\documentclass[12pt,reqno]{amsart}
\usepackage[english]{babel}
\usepackage[cp1251]{inputenc}

\usepackage{mathrsfs}
\usepackage{amsthm}
\usepackage{amssymb}
\usepackage{amsmath,amsthm,amsfonts}
\usepackage{mathabx}
\usepackage{hhline}
\usepackage{textcomp}
\usepackage{amsmath,amscd}
\usepackage[all,cmtip]{xy}
\usepackage{mathtools}
\usepackage[margin=1in]{geometry}
\usepackage{lipsum}
\usepackage{stmaryrd}
\usepackage{extarrows}

\author[]{Kaveh Eftekharinasab}
\title[]{Fr\'{e}chet Lie algebroids and their cohomology}
\address{Topology dept. \\ Institute of Mathematics of NAS of Ukraine \\ Te\-re\-shchen\-kivska st. 3, Kyiv, 01601 Ukraine}
\email{kaveh@imath.kiev.ua}
\makeatletter
\@addtoreset{equation}{section}
\@addtoreset{figure}{section}
\@addtoreset{table}{section}
\makeatother

\newtheorem{theorem}{Theorem}[section]

\newtheorem{defn}{Definition}[section]

\theoremstyle{definition}

\DeclareMathAlphabet{\mathpzc}{OT1}{pzc}{m}{it}

\newcommand{\rr}{\mathbb{R}}
\newcommand{\nn}{\mathbb{N}}
\newcommand{\bb}{\mathcal{B}}
\newcommand{\xx}{\mathfrak{X}}
\newcommand{\lx}{\mathcal{L}}
\newcommand{\sn}{\mathrm{p}_n}
\newcommand{\sno}{\mathrm{P}_B^n}
\newcommand{\td}{\mathbb{T}}
\DeclareMathOperator{\dd}{d}
\DeclareMathOperator{\im}{im}
\newcommand\bGroup{\mathbf{H}}

\begin{document}

\begin{abstract}
We define Lie and Courant algebroids on Fr\'{e}chet manifolds. Moreover, we construct
a Dirac structure on the generalized tangent bundle of a Fr\'{e}chet manifold and show that it 
inherits a Fr\'{e}chet Lie algebroid structure. We show that the Lie
algebroid cohomology of the $\bb$-cotangent bundle Lie algebroid of a
weakly symplectic Fr\'{e}chet manifold $M$ is the Lichnerowicz-Poisson cohomology of $M$.

\end{abstract}
\maketitle
\section{Introduction}
In recent years Poisson geometry has been extended to the Banach manifolds context. In particular,
the concept of Lie algebroid was generalized to the category of Banach vector bundles in~\cite{ans,f}. 
Lie algebroids were also defined on Projective limits of Banach Manifolds~\cite{ca}. Dirac structures on Banach
manifolds were studied in~\cite{an2,v}. 
The assertion that the Lie algebroid cohomology of the cotangent bundle Lie algebroid of a finite dimensional Poisson manifold $M$  
is the Lichnerowicz-Poisson cohomology of $M$  was generalized to the Banach manifolds case in~\cite{i}.

Our goal in this paper is to extend to Fr\'{e}chet manifolds some of the aforementioned results. 
Frechet manifolds arise in number of problems that have significance in global analysis and physical field theory. 
However, due to permanent problems with Fr\'{e}chet spaces (i.e problems of intrinsic nature) in most cases they are 
handled by indirect methods or only certain type of Fr\'{e}chet manifolds are considered 
(see~\cite{dd} for a survey on recent developments in Fr\'{e}chet geometry).  
One of the main issues in the theory of Fr\'{e}chet spaces is that the dual of a proper Fr\'{e}chet space (not Banachable) is
never a Fr\'{e}chet space. In addition, the space of continuous linear mappings of one Fr\'{e}chet space
to another is not a Fr\'{e}chet space in general. This defect puts in question the way of defining cotangent bundle.
In fact, as pointed out in~\cite{neeb}, if $M$ is a manifold modelled on a Fr\'{e}chet space $F$, then
there is no vector topology on the dual of $F$ that can lead to a smooth manifold structure on the set-theoretic
cotangent bundle of $M$.
But the cotangent bundle of a Poisson manifold is a special case of Lie algebroid and vital to the Cartan calculus
of differential forms. A way out of this difficulty was observed in~\cite{tw},
where a notion of the $\bb$-cotangent bundle of a manifold modelled on a locally convex space was introduced
and the Cartan calculus of differential forms was successfully adapted to the category of manifolds modelled on locally
convex spaces. 

By means of this notion we define Lie and Courant algebroids on Fr\'{e}chet manifolds. 
Furthermore, we construct a Dirac structure as a subbundle of the generalized tangent bundle $TM_{\bb}'\oplus TM$ of a Fr\'{e}chet manifold 
$M$ and show 
that it inherits a Lie algebroid structure from the Courant bracket. We show that a weak symplectic form of a Fr\'{e}chet manifold
$M$ determines the so-called Lichnerowicz-Poisson cohomology of $M$ and the Chevalley-Eilenberg cohomology of its $\bb$-cotangent bundle
Lie algebroid which are exactly the same.

We should mention that another approach to  the geometry of Fr\'{e}chet cotangent bundle is using the convenient setting
which provides two notions of cotangent bundle (kinematic and operational).
We can attempt to use convenient calculus to develop Poisson geometry for Fr\'{e}chet manifolds but in 
this paper we consider only Micheal-Bastiani differentiability which is more familiar and applicable 
for people working in Fr\'{e}chet spaces. It turns out by using the notion of $\bb$-cotangent bundle most of the assertions and constructions are much the same as those of
Banach manifolds case.

\section{Preliminaries: Poisson Fr\'{e}chet manifolds}

In this section, following~\cite{tw}, we define Poisson structures on Fr\'{e}chet manifolds.
We will apply the notion of differentiability in the Michal-Bastiani sense. We will be working in the category of
smooth manifolds and bundles.
\begin{defn}\label{di}
 Let $E$ and $F$ be Fr\'{e}chet spaces, $U \subseteq E$  open and $f : U \rightarrow F$ a continuous map.
 The derivative of $f$ at $x \in U$ in the direction of $h \in E$ is defined as
 $$
 \dd f(x)(h) \coloneq \lim_{t \rightarrow 0} \frac{1}{t} (f(x+ht) - f(x))
 $$
 whenever the limit exits. The map $f$ is called differentiable at $x $ if $\dd f(x)(h)$ exists for all $h \in E$. It is called a
 $C^1$-map if it is differentiable at all points of $ U$ and 
 $$
 \dd f : U \times E \rightarrow F, \quad (x,h) \mapsto \dd f(x)(h)
 $$
 is a continuous map. Higher directional derivatives and $C^k$-maps, $k\geq 2$, are defined in the obvious inductive fashion. 
\end{defn}
Within this framework Fr\'{e}chet manifolds, Fr\'{e}chet vector bundles (especially tangent bundles) and $C^k$-maps between
Fr\'{e}chet manifolds are defined in obvious way (cf.~\cite{neeb}).
However, for a manifold $M$ modelled on a Fr\'{e}chet space $F$ we can define the set-theoretic cotangent bundle $T'M$ 
(without any topology on the fiber), but in general there
is no vector topology on $F'$, the dual of $F$, that can lead to the identification $T'M \cong F \times F'$, see~\cite[Remark I.3.9]{neeb}. 
Thus, we follow~\cite{tw} and use the notion of a $\bb$-cotangent bundle instead. In this definition 
to put a manifold structure on $T'M$, the dual of $F$ is equipped by a $\bb$-topology, where $\bb$ is a
bornology on $F$. To be precise, we recall that a family $\bb$ of bounded subsets of $F$ that covers $F$
is called a bornology on $F$ if it is directed upwards by inclusion and if for every $B \in \bb$ and $r \in \rr$ there is
a $C \in \bb$ such that $r \cdot B \subset C$. 

Let $E$ be a Fr\'{e}chet space, $\bb$ a bornology on $E$ and $L_{\bb}(E,F)$  the space of all linear continuous maps from $E$ to $F$. 
The $\bb$-topology on $L_{\bb}(E,F)$ is a Hausdorff locally convex topology defined by all seminorms 
 $ \sno (L) \coloneq \sup \{ \sn (L (e)) \mid  e \in B\}$, 
where $B \in \bb$ and $\{ \sn \}_{ n\in \nn}$ is a family of seminorms defining the topology of $F$. 
One similarly may define $L_{\bb}^k(E,F)$ and $\bigwedge^{\mathclap{k}} L _{\bb}(E,\rr)$, the space of $k$-linear jointly continuous  maps of
$E^k$ to $F$ and the space of anti-symmetric $k$-linear jointly continuous  maps of $E^k$  to $\rr$, respectively. If $\bb$ contains all compact sets, 
then the $\bb$-topology on the space  $L _{\bb}(E,\rr) = E_{\bb}'$ of all continuous linear functional on $E$, the dual
of $E$, is the topology of compact convergence. 

If $\bb$ contains all compact sets of $E$, then we define the differentiability of class $C^k_{\bb}$ : Let $U \subset E$ be open, a map 
$f : U \rightarrow F$ is called $C_{\bb}^1$ if its partial derivatives exist and the induced map
$\dd f : U \rightarrow L_{\bb}(E,F)$ is continuous. Similarly
we can define maps of class $C^k_{\bb}$, $k \in \nn \cup \{\infty\}$, see~\cite[Definition 2.5.0]{ke}. 
A map $f : U \rightarrow F$ is $C^{k}_{\bb} \,,k\geqq 1,$ if and only if  $f$ is $C^k$ in the sense of
Definition~\eqref{di}, see~\cite[Theorem 2.7.0 and Corollary 1.0.4 (2)]{ke}.
In particular, $f$ is $C^{\infty}_{\bb}$ if and only in $f$ is $C^{\infty}$. Thus, if $f$ at $x \in E$ is $C^k$ and hence $C^k_{\bb}$,
the derivative $f$ at $x$, $\dd f(x)$, is an element of $E_{\bb}'$.  

Assume that $\bb$ is a bornology on $F$ containing all compact sets and $M$ is a Fr\'{e}chet manifold  
modelled on $F$. Let $f$ be a functional defined over $M$. The derivative of $f$ at $x \in M$ can be written in  terms of the iterated tangent bundles of $M$ and we can consider 
$\dd f: TM \rightarrow F$ given by $\dd f(x,h) = \dd f(x)(h)$ upon locally identifying $TM$
with $U \times F$, where $U$ is an open set in $F$. Therefore, if at $x \in M$ a map $f : M \rightarrow \rr$ is $C^{k}$ and hence 
$C_{\bb}^{k}$, then $\dd f(x)$ belongs to $L_{\bb}(T_xM,\rr)=(T_xM)'_{\bb}$.
\begin{defn}
Let $M$ be a Fr\'{e}chet manifold modelled on a Fr\'{e}chet space $F$ and $\bb$ a bornology on $F$.
 The $\bb$-cotangent bundle of $M$ is defined as $TM_{\bb}' \coloneq \bigcup_{x \in M}(T_xM)'_{\bb} $ and the
 $k$-exterior product of the $\bb$-cotangent bundle as
 $\bigwedge^{\mathclap{k}} TM_{\bb}' \coloneq \bigcup_{x \in M}\bigwedge^{\mathclap{k}} (T_xM)'_{\bb}$.
\end{defn}
If $\bb$ is chosen such that $T(\bb) \subset \bb$ for all continuous linear endomorphisms $T$ of $F$, then $\bigwedge^{\mathclap{k}} TM_{\bb}' $ 
is a vector bundle in the category of locally convex spaces with the local model $F \times \bigwedge^{\mathclap{k}}F_{\bb}'$. In particular,
$TM_{\bb}'$ is a vector bundle in the category of locally convex spaces with the  local model $F \times F_{\bb}'$, see~\cite[Remark (1), p. 339]{tw}.
Therefore, we always assume that bornologies have this property and contain all compact sets.

Let $M$ be a manifold modelled on a Fr\'{e}chet space $F$ and $\bb$ a bornology on $F$.
A smooth differentiable $k$-form of type $\bb$ is a smooth section of the bundle $\bigwedge^{\mathclap{k}} TM_{\bb}'$.
We can also define a smooth differential $k$-form in the weak sense (which is usually used in the literature) as a section of the set-theoretic $k$-exterior bundle $\bigwedge^{\mathclap{k}}TM'$, 
cf.~\cite{neeb}. A section $\omega$ of $\bigwedge^{\mathclap{k}}TM_{\bb}' \rightarrow M$ is a smooth differential $k$-form of type $\bb$ if and only if 
$\omega$ is a smooth differential $k$-form in the weak sense~\cite[Proposition IV.6]{tw}. Thus, in the sequel we call smooth
differential forms of type $\bb$ simply smooth differential forms. We always assume that differential
forms are smooth without mentioning it and write $\omega_x$ instead of $\omega(x)$ for $x \in M$. 

We denote by $\xx(M)$ the space of all vector fields on $M$. The Lie bracket $[X,Y]$ of $X,Y \in \xx(M)$ 
 is again vector field and $(\xx(M),[\cdot,\cdot])$ is a Lie algebra~\cite[Proposition II.3.7]{neeb}.
 We further obtain for each $X \in \xx(M)$ and a $k$-form $\omega$ on $M$ a unique linear map
 $ (\imath_X \omega)_{x} = \imath_{X(x)}\omega_x$, where $x \in M$ and $v,v_0, \cdots ,v_{k-1} \in T_xM$ and $(\imath_v \omega_x)(v_1,\cdots,v_{k-1}) \coloneq \omega_x(v,v_1,\cdots,v_{k-1})$.
 Let $\omega$ be a  $k$-form on $M$. Then $\dd_{\mathrm{dR}} \omega$, the de Rham derivative of $\omega$, on vector fields
 $X_0,\cdots,X_k \in \xx(M)$ is a smooth
 $(k+1)$-form~\cite[Lemma IV.8]{tw} and is given by
 \begin{align} \label{dr}
 (\dd_{\mathrm{dR}}\omega)(X_0, \cdots,X_k) &= \sum_{i=0}^{i=k}(-1)^i X_i(\omega(X_0,\cdots,\hat{X_i},\cdots,X_k))+ \nonumber\\
 &+ \sum_{0\leqq i < j \leqq k} (-1)^{(i+j)} \omega([X_i,X_j],X_0,\cdots,\hat{X_i},\cdots, \hat{X_j}, \cdots, X_k),
\end{align}
where a hat over symbols means omission. We now define the Lie derivative $\lx_X$ of a differential form in the direction
of a vector field $X$ by the Cartan formula:
\begin{equation}
 \lx_X =\dd_{\mathrm{dR}} \circ \, \imath_X + \imath_X \circ \dd_{\mathrm{dR}}.
\end{equation}
\begin{defn}
 A Fr\'{e}chet manifold $M$ is called weakly symplectic if for a closed smooth 2-form $\omega$ on $M$ the linear continuous map 
 \begin{equation*}
\begin{array}{cccc}
\omega^{\hash}:T_xM \longrightarrow (T_xM)_b' \\
v_x \mapsto \omega_x(v_x,\cdot) =\imath_{v_x} \omega_x
\end{array}
\end{equation*}
is injective for all $x\in M$. Here, $(T_xM)_b'$ is the strong dual of the tangent space.
\end{defn}
\begin{defn}
 A vector field $X_f$ on a weakly symplectic Fr\'{e}chet manifold $(M,\omega)$ is called
 the symplectic gradient vector field of a smooth real valued function $f \in C^{\infty}(M)$ if 
 $\dd f = -\omega(X_f,\cdot) = \omega^{\hash}(-X_f) $. Let $f,g \in C^{\infty}(M)$ be such that $X_f$ and $X_g$ exist. 
 For $f$ and $g$, the Poisson structure  $\{f,g\}$ of is defined by
 \begin{equation}
  \{f,g\} \coloneq \omega (X_f,X_g).
 \end{equation}
 It is $\rr$-linear, anti-symmetric and satisfies Jacobi identity, if all involved symplectic gradient vector fields exist.
 We say that a pair $\{M,\{\cdot,\cdot\}\}$ is a Fr\'{e}chet Poisson manifold.
\end{defn}

\section{Fr\'{e}chet Lie  algebroids}

Let $M$ be a manifold modelled on a Fr\'{e}chet space $F$ and let $\pi : L \rightarrow M$ be a Fr\'{e}chet vector bundle over $M$ with fibers of type
$\mathbb{F}$. We denote by $\Gamma(L)$ the space of smooth sections of the vector bundle $L$. The spaces $\Gamma(L)$ and $\xx(M)$ are both
$C^{\infty}(M)$-modules. 
\begin{defn}\label{th:cl} A Lie algebroid $L$ over $M$ is a vector bundle $\pi : L \rightarrow M$ together with 
a bracket $[\cdot,\cdot]_L$ on the space $\Gamma (L)$  and a bundle map $\lambda_L: L \rightarrow TL$, called anchor, such that 
 
\begin{enumerate}
 \item The induced map $\lambda_L : (\Gamma(L),[\cdot,\cdot]_L ) \rightarrow (\xx(M),[\cdot,\cdot])$ 
 given by $(\lambda_L (s))(x) = \lambda_L (s(x))$, $ x\in M$, $s \in \Gamma(L)$ is a Lie algebra homomorphism,
 \item $[s_1,fs_2]_L = f [s_1,s_2]_L + \lambda_L(s_1)(f)s_2$ for every $f \in C^{\infty}(M)$ and $s_1,s_2 \in \Gamma(L)$.
 \end{enumerate}
 \end{defn} 
 The tangent bundle $TM$ is trivially a Fr\'{e}chet Lie algebroid for the usual Lie
 bracket of vector fields on $M$ and the identity map of $TM$ as an anchor map.

 \begin{defn}
  A Courant algebroid is a vector bundle $\pi :  \mathsf{C} \rightarrow M $ together with an anchor $\lambda_{\mathsf{C}}$,
  a nondegenerate symmetric bilinear form $\Theta$ and a bracket $[\cdot,\cdot]_{\mathsf{C}}$ on $\Gamma(\mathsf{C})$ such that
  for all $s_1,s_2,s_3 \in \Gamma(\mathsf{C})$ and $f \in C^{\infty}(M)$
  \begin{enumerate}
   \item $[s_1,[s_2,s_3]_{\mathsf{C}}]_{\mathsf{C}} = [[s_1,s_2]_{\mathsf{C}},s_3]_{\mathsf{C}} + [s_2,[s_1,s_3]_{\mathsf{C}}]_{\mathsf{C}}$,
   \item $ \lambda_{\mathsf{C}}(s_1)\Theta(s_2,s_3) = \Theta([s_1,s_2]_{\mathsf{C}},s_3) + \Theta[s_2,[s_1,s_3]_{\mathsf{C}}]_{\mathsf{C}}$,
   \item $ [s_1,s_2]_{\mathsf{C}} + [s_2,s_1]_{\mathsf{C}} = \Delta (\Theta (s_1,s_2))$, where $\Delta: C^{\infty}(M) \rightarrow \Gamma(\mathsf{C})$
   is defined by $\Theta(\Delta(f), s)$ $ = \lambda_{\mathsf{C}}(s)f$.
  \end{enumerate}
 \end{defn}
Since the $\bb$-cotangent bundle $TM_{\bb}'$ of $M$ is a vector bundle in the category of locally convex spaces with the local model
$F \times F_{\bb}'$, then the Whitney sum $TM \oplus TM_{\bb}'$ makes sense. We denote by $\mathbb{T}M = TM \oplus TM_{\bb}' $ the generalized
tangent bundle of $M$. Let $X,Y \in \xx(M)$ and $\alpha, \beta \in \bigwedge^{\mathclap{1}} TM_{\bb}' $. Now define
the bracket $[X,Y]_{\mathbb{T}M} = ([X,Y],\lx_X \beta - (\imath_Y \circ \dd_{\mathrm{dR}} \alpha))$ and the anchor $\lambda_{\mathbb{T}M}$ given by
$\lambda_{\mathbb{T}M}(X,\alpha) = X$. If we define $\Delta_{\mathbb{T}M} ((X,\alpha), (Y,\beta)) = \alpha(Y) + \beta(X)$, then
we can easily verify that for $(\mathbb{T}M,[X,Y]_{\mathbb{T}M},\lambda_{\mathbb{T}M})$ the conditions (1)-(2) of the Definition~\ref{th:cl} are fulfilled, therefore, $\mathbb{T}M$ is a Courant algebroid.
The orthogonal complement $\mathbb{L}^{\bot}$ of the subbundle $\mathbb{L} \subset \mathbb{T}M$ is defined as follows
$$
\mathbb{L}^{\bot} \coloneq \left\lbrace (X,\alpha) \in \mathbb{T}M : \Delta_{\mathbb{T}M} ((X,\alpha), (Y,\beta)) =0, \forall (Y,\beta) \in \mathbb{T}M \right\rbrace
.$$
\begin{defn}
 A vector subbundle $\mathbb{D}$ of the Courant algebroid $\mathbb{T}M$
 that coincides with its  orthogonal complement $\mathbb{D}^{\bot}$  with respect to $\Delta_{\mathbb{TM}}$ is said to be
 an almost Dirac structure. It is called a Dirac structure if, in addition, is closed under the bracket $[\cdot,\cdot]_{\mathbb{TM}}$.
\end{defn}
Define the Courant bracket on $\Gamma(\td M)$ by
$$
\llbracket (X,\alpha), (Y,\beta) \rrbracket = \left( [X,Y], \lx_X \beta - \lx_Y \alpha + \dfrac{1}{2} \dd_{\mathrm{dR}} (\alpha(Y) - \beta(X)) \right).
$$
We can easily show that the restriction of $  \llbracket \cdot,\cdot \rrbracket$ to $\Gamma(\mathbb{D})$ yields a Lie bracket and if
we let $\mathrm{Pr} : \mathbb{D} \rightarrow TM$ to be the restriction of the projection to $TM$, then $(\mathbb{D},\llbracket \cdot,\cdot \rrbracket_{\mathbb{D}},\mathrm{Pr}) $
is a Fr\'{e}chet Lie algebroid.

\section{Fr\'{e}chet Lie algebroids cohomology}

Let $(M,\omega)$ be a weakly symplectic Fr\'{e}chet manifold. We denote by $\xx^k(M)$ and $\Omega^k(M)$ the spaces of
all $k$-vector fields and $k$-differential forms on $M$, respectively. Define a morphism
\begin{equation}\label{anch}
 \hash \omega : \Omega^1(M) \rightarrow \xx^1(M); \, \beta (\hash \omega (\alpha)) = \omega (\alpha,\beta), \, \forall \alpha, \beta \in \Omega^1(M). 
\end{equation}
The weak symplectic form $\omega$ induces a unique lie bracket of $1$-forms given by
\begin{equation}\label{bra}
 \{\alpha,\beta\} = \lx_{\hash \omega(\alpha)}\beta - \lx_{\hash \omega(\beta)}\alpha - \dd_{\mathrm{dR}}\omega(\alpha,\beta).
\end{equation}
In general the existence of the lie bracket is equivalent to the existence of the  weak symplectic form $\omega$,
the proof is the same as the finite dimensional case, see~\cite{bh}.
 Define the contravariant exterior differential $ \sigma: \xx^k(M) \rightarrow \xx^{k+1}(M)$ for $X \in \xx^k(M)$ by
 \begin{align}\label{sigma}
 (\sigma X)(\alpha_0, \cdots,\alpha_k) &= \sum_{i=0}^{i=k}(-1)^i \hash \omega (\alpha_i)(X(\alpha_0,\cdots,\hat{\alpha_i},\cdots,\alpha_k))+ \nonumber\\
 &+ \sum_{0\leqq i < j \leqq k} (-1)^{(i+j)} X(\{\alpha_i,\alpha_j\},\alpha_0,\cdots,\hat{\alpha_i},\cdots, \hat{\alpha_j}, \cdots, \alpha_k),
\end{align}
where $\alpha_0 \cdots, \alpha_k \in \Omega(M)$ and a hat over symbols means omission. 
Formally, the expression~\eqref{sigma} is exactly the same as the de Rham derivative of forms and hence
its algebraic consequences will be the same, in particular
\begin{description}
 \item [(i)] $\sigma^2 =0$, 
 \item[(ii)] $\sigma (X_1 \wedge X_2) = \sigma(X_1) \wedge X_2 + (-1)^{\deg X_1}X_1 \wedge \sigma(X_2)$,
 \item[(iii)] $\sigma ([X_1 , X_2]) = -[ \sigma(X_1), X_2] - (-1)^{\deg X_1}[X_1, \sigma(X_2)]$,
\end{description}
where $X_i (i=1,2)$ are $k$-vector fields on $M$ and $\deg X^i$ is the degree of $X^i$.
Therefore, $ \Omega (M) \coloneq \displaystyle \oplus_{k\in \nn \cup \{0\}}(\Omega^k M) $ with the coboundary operator $\sigma$
is cochain complex and we can define the Lichnerowicz-Poisson cohomology of $M$.
\begin{defn}\label{co}
 $\Omega (M)$ with the coboundary operator $\sigma$ is called the Lichnerowicz-Poisson cochain of $M$, and
 \begin{equation}
  \bGroup^k_{LP} (M,\omega) \coloneq \frac{\ker \left(\xx^{k}(M) \xlongrightarrow{\sigma} \xx^{k+1}(M) \right)}
  {\im \left(\xx^{k-1}(M) \xlongrightarrow{\sigma} \xx^{k}(M) \right)}
 \end{equation}
are Lichnerowicz-Poisson cohomology or LP-cohomology spaces of $M$.
\end{defn}
The following LP-cohomology spaces can be computed straightaway by Definition~\ref{co}.

Let 
$ Z_{ \{C^{\infty}(M),\{,\} \}} = \left\lbrace f \in C^{\infty}(M) : \forall g \in C^{\infty}(M),\, X_gf=0 \right\rbrace$
and let $\xx^1_{S}(M)$ be the space of symplectic gradient vector fields $X_f,\, f \in C^{\infty}(M)$. 
Let $\xx^1_{\omega}(M) \coloneq \left\lbrace X \in \xx(M): \lx_X \omega =0   \right\rbrace$.
We then have
\begin{itemize}
 \item [(i)] $\bGroup^0_{LP}(M,\omega) = Z_{ \{C^{\infty}(M),\{,\} \}}$, since $\sigma f =-X_f$. 
 \item [(ii)] $ \bGroup^1_{LP}(M,\omega) = \dfrac{\xx^1_{\omega}(M)}{\xx^1_{S}(M)}$.
\end{itemize}
Now we define the Chevally-Eilenberg cohomology~\cite{ce} associated to Fr\'{e}chet Lie algebroids. 
Let $(L,\lambda_L,[\cdot,\cdot]_E)$ be a Fr\'{e}chet Lie algebroid over a Fr\'{e}chet manifold $M$. 
Let $C^{\infty}(M)$ act on $\Gamma(L)$ by $(s,f) \mapsto \lambda_L(s)f$.
A $k$-linear anti-symmetric mapping $\ell_k :\Gamma(L)^k \rightarrow C^{\infty}(M)$ is called a $C^{\infty}(M)$-valued $k$-cochain.
Let $C^k \left( \Gamma(L);C^{\infty}(M) \right)$ be the vector space of these cochains.
Define the operator $\dd_L$
by
\begin{align*} \
 (\dd_L \ell_k)(s_0, \cdots,s_k) &= \sum_{i=0}^{i=k}(-1)^i \lambda_L(s_i)(\ell_k (s_0,\cdots,\hat{s_i},\cdots,s_k))+\\
 &+ \sum_{0\leqq i < j \leqq k} (-1)^{(i+j)} \ell_k ([s_i,s_j]_L,s_0,\cdots,\hat{s_i},\cdots, \hat{s_j}, \cdots, s_k),
\end{align*}
for a $k$-cochain $\ell_k$ and $s_0, \cdots, s_k \in \Gamma(L)$. Like the case of the de Raham derivative of forms we obtain $\dd_L \circ \dd_L = 0$.
Therefore, $C^k \left( \Gamma(L);C^{\infty}(M) \right)$ with $\dd_L$ forms a Chevalley-Eilenberg cochain and
the corresponding cohomology spaces $$
\bGroup^k (C^k \left( \Gamma(L);C^{\infty}(M) \right)) \coloneq \frac{\ker \left(\Gamma(L)^{k} \xlongrightarrow{\dd_L}
\Gamma(L)^{k+1} \right)}{\im \left(\Gamma(L)^{k-1} \xlongrightarrow{\dd_L} \Gamma(L)^{k} \right)}
,$$
are called the Lie algebroid cohomology of $\Gamma(L)$ with coefficient in $C^{\infty}(M)$. 
For the tangent bundle Lie algebroid $TM$ of $M$,  the Lie algebroid cohomology 
is just the de Raham cohomology of $M$.

On any weakly symplectic Fr\'{e}chet manifold $(M,\omega) $ the bracket $\{\cdot,\cdot\}$ of $1$-forms~\eqref{bra} defines a Lie algebroid
structure on the $\bb$-cotangent bundle $TM_{\bb}'$ with the anchor $\hash \omega: TM_{\bb}' \rightarrow TM $ given
by $\beta(\hash \omega(\alpha)) = \omega(\beta,\alpha); \, \alpha, \beta \in TM_{\bb}'$. In this case, eventually 
the operator
$\dd_{(TM_{\bb}')}$ coincides with the contravariant exterior differential~\eqref{sigma} and so we have
\begin{theorem}
 The Lie algebroid cohomology of the $\bb$-cotangent bundle Lie algebroid is the Lichnerowicz-Poisson cohomology
 of $M$.
\end{theorem}

This result  was obtained  for finite dimensional manifolds in~\cite{bh}.

\end{document}